\documentclass{amsart}
\usepackage{amssymb}

\usepackage{tikz-cd}

\let\ge\geqslant

\def\1{^{-1}}

\def\CP{{\mathbf C\mathbf P}}
\def\MU{{\mathbf M\mathbf U}}
\def\MSU{{\mathbf M\mathbf S \mathbf U}}

\newtheorem{theorem}{Theorem}[section]
\newtheorem{lemma}[theorem]{Lemma}
\newtheorem{proposition}[theorem]{Proposition}

\theoremstyle{definition}

\newtheorem{corollary}[theorem]{Corollary}

\theoremstyle{remark}
\newtheorem{remark}[theorem]{Remark}

\newcommand{\BU}{\mathop{\mathrm{BU}}}

\numberwithin{equation}{section}

\begin{document}



\commby{}


	
\title[]{Polynomial generators of $\MSU^*[1/2]$ related to classifying maps of certain formal group laws}

\begin{abstract}
This paper presents a commutative complex oriented cohomology theory that realizes the Buchstaber formal group law $F_B$ localized away from 2. It is shown that the restriction of the classifying map of $F_B$ on special unitary cobordism ring localized away from 2 defines a four parameter genus, studied by Hoehn and Totaro.    
\end{abstract}

\author{Malkhaz Bakuradze }
\address{Faculty of exact and natural sciences, A. Razmadze Math. Institute, Iv. Javakhishvili Tbilisi State University, Georgia }
\email{malkhaz.bakuradze@tsu.ge } 
\thanks{The author was supported by Shota Rustaveli NSF grant FR-21-4713}
\subjclass[2010]{55N22; 55N35}
\keywords{Complex bordism, $SU$-bordism,  Formal group law, Complex elliptic genus}

\maketitle

\section{Introduction}

The ring of complex cobordism $\MU_*$ and  the ring $\MSU_*$ of special unitary cobordism has been studied by many authors.
We refer the reader to \cite{Stong}, \cite{L-C-P} for details. In particular the ring $\MSU_*$, localized away from 2, is torsion free
$$\MSU_*[1/2]=\mathbb{Z}[1/2][x_2,x_3,\cdots],\,\,\,\,|x_i|=2i$$
and $SU$-structure forgetful homomorphism is the inclusion in complex cobordism ring
$$\MSU_*[1/2]\subset \MU_*[1/2]=\mathbb{Z}[1/2][x_1,x_2,x_3,\cdots].$$  

\medskip

In this paper we construct a commutative complex oriented cohomology theory (Theorem \ref{main-lemma}) such that the coefficient ring is the scalar ring of the Buchstaber formal group  law $F_B$ with inverted 2, and show (Proposition \ref{main-lemma})
that after restricted to $\MSU_*[1/2]$, the classifying map of $F_B$ can become a genus 
\begin{equation}
	\label{eq:1-1}
	\MSU_*[1/2]\to \mathbb{Z}[1/2][x_2,x_3,x_4],\,\,\,\,|x_i|=2i,
\end{equation}
studied by Hoehn \cite{H} and Totaro  \cite{TO}. 

Since $\MSU^*$ is not complex oriented, it is difficult to compute the
genus \eqref{eq:1-1} on specific explicit elements. Using the polynomial generators of the spherical cobordism
ring $W^*[1/2]$ given by G. Chernykh and T. Panov \cite{C-P}, we derive certain polynomial generators of $\MSU^*[1/2]$ in terms of the universal formal group law. This gives a new understanding of the genus \eqref{eq:1-1}.

\medskip

In particular, the classifying map $f_B$ of $F_B$ is a surjection on some infinitely generated ring $\Lambda_B$, with kernel generated by some explicit elements 
(Proposition \ref{prop.1}).  
After it is tensored with rationals it is identical (Proposition \ref{prop.2}) to the complex elliptic genus
\begin{equation}
\label{1-2}
\MU_*\otimes \mathbb{Q}\to \mathbb{Q}[x_1,x_2,x_3,x_4],
\end{equation}
defined in Hoehn’s thesis \cite{H}(Section 2.5).
Here $x_1$ is the image of complex projective plane $\CP_1$ and
$x_2, x_3, x_4$ are the images of any first three generators of the polynomial ring $\MSU_*\otimes \mathbb{Q}=\mathbb{Q}[x_2,x_3,x_4,\cdots]$.

By Hoehn \cite{H}, for $X$ an
$SU$-manifold of complex dimension $n$, the exponential characteristic class $\phi(X)$ is in fact a Jacobi form of weight
$n$. Jacobi forms are generalizations of modular forms. See details in \cite{TO}.
Hoehn showed that the Jacobi forms $x_2, x_3, x_4$ arise as the elliptic
genera of certain explicit $SU$-manifolds, of complex dimensions 2, 3, 4, so that 
the homomorphism
$$\MSU_* \to (\text{ Jacobi forms over } \mathbb{Z})$$
becomes surjective after it is tensored with $\mathbb{Z}[1/2]$.

In \cite{TO} (Theorem 4.1) Totaro proved that the Krichever-Hoehn complex elliptic genus on complex cobordism  viewed as a homomorphism \eqref{1-2}
is surjective and the kernel is equal to the ideal of complex flops. 

Then Totaro proved (Theorem 6.1) that the kernel of the complex elliptic genus on $\MSU_*\otimes \mathbb{Z}[1/2]
$ is equal to the ideal $I$ of $SU$-flops. Also, the quotient
ring is a polynomial ring:
\begin{equation}
\label{genus-TO}
\MSU_*[1/2]/I = \mathbb{Z}[1/2][x_2, x_3, x_4].
\end{equation}

Unfortunately  $\MSU_*$ is not complex oriented.  It would be nice to develop a method for calculating \eqref{genus-TO} explicitly in terms of the universal formal group law using some generators of $\MSU_* [1/2]$ treated as explicit elements in $\MU_*[1/2 ]$.

This goal can be achieved as follows: in Section \ref{5} we replace the ideal of complex flops with a more explicit ideal by considering the integral Buchstaber genus which is identical to the Krichever-Hoehn complex elliptic genus over $\MU_*\otimes \mathbb{Q}$.

In Section \ref{6} we use the polynomial generators of the spherical cobordism ring $W_*[1/2]$ constructed in \cite{C-P} and define certain polynomial generators of $\MSU_*[1/2]$. 
In particular, we use fact that $W_*[1/2] $ is generated by the coefficients of the corresponding formal group law.  Given Novikov's criteria and that $W_*[1/2]$ is a free $\MSU_*[1/2]$ module generated by $1$ and $\CP^1$, we define some generators in $\MSU_*[1/2]$ ( Proposition \ref{msu-b_k}) in terms of the universal formal group law. Finally the explicit quotient map of the Buchstaber formal group law (Proposition \ref{prop.1}) gives the decompositions of constructed generators $\MSU_*$ of dimensions $>10$ in polynomial ring $\mathbb{Z}[1/2][x_2,x_3,x_4]$. It is a way to calculate the genus \eqref{genus-TO} in terms of the universal formal group law.

In Section \ref{7} we consider the restriction of classifying map of the universal abelian formal group law on $\MSU_*[1/2]$ to define a genus with one parameter.

\bigskip

\section{Preliminaries}
\label{2}

 The theory $W_*$ of $c_1$-spherical bordism is defined geometrically in \cite{Stong} (Chapter VIII). The closed manifolds $M$ with a $c_1$-spherical structure, consist of

– a stably complex structure on the tangent bundle $T M$;

– a $\CP^1$-reduction of the determinant bundle, that is, a map $f : M \to \CP^1$
and an
equivalence $f^*(\eta)\simeq \det T M$, where $\eta$ is the tautological bundle over $\CP^1$.

This is a natural generalization of an $SU$-structure, which can be thought of as a trivialization of the determinant bundle. The corresponding bordism theory is called $c_1$
-spherical bordism and is denoted $W_*$.
The unitary and special unitary bordism rings are denoted by $\MU_*$ and $\MSU_*$ respectively. We refer to \cite{L-C-P} and \cite{C-P}  for details on $\MSU_*$ and $W_*$ that will be used throughout the paper. 

Motivated by string theory in  \cite{KR}, \cite{H}, \cite {TO} the universal Krichever-Hoehn complex elliptic genus $\phi_{KH}$ is defined as the ring homomorphism 
\begin{equation}
\label{0}
\phi_{KH}:\MU^* \to \mathbb{Q}[q_1,q_2,q_3,q_4]
\end{equation}
associated to the  Hirzebruch characteristic power series $Q(x)=\frac{x}{f(x)}$,
where $$h(x)=\frac{f'(x)}{f(x)}$$
is the solution of the differential equation in $\MU^*\otimes \mathbb{Q}$
\begin{equation}
\label{Hoehn}
(h')^2=S(h),
\end{equation}
where
$$
S(x)=x^4+q_1x^3+q_2x^2+q_3x+q_4,
$$ 
for some formal parameters $|q_i|=2i$.

\medskip

One consequence of Krichever-Hoehn’s rigidity theorem  \cite{KR}, \cite{H},  \cite{TO},  is that (\cite{H}, Kor 2.2.3) if $F \to E \to B$ is a fiber
bundle of closed connected weakly complex manifolds, with structure group a
compact connected Lie group $G$, and if $F$ is an $SU$-manifold, then the elliptic
genus $\phi_{KH}$ satisfies $\phi_{KH}(E) = \phi_{KH}(F)\phi_{KH}(B)$. In fact, the elliptic genus is the universal genus with the above multiplicative property

In  (\cite{TO}, Theorem 6.1) Totaro gave a geometric description of the kernel ideal $I$ of the complex
elliptic genus restricted to $\MSU[1/2]$, the ideal of $SU$ flops.  This kernel is equal to the ideal in $\MSU_*[1/2]$ generated by twisted projective bundles $\tilde{\CP}(A \oplus B)$
over weakly complex manifolds $Z$ such that the complex vector bundles $A$
and $B$ over $Z$ have rank $2$ and $c_1Z + c_1A + c_1B=0$; in this case, the total
space is an $SU$-manifold. Then Totaros's result says that the
$I\in \MU_*[1/2]$ contains a polynomial generator of $\MSU_*[1/2]$ in real dimension $2n$
for all $n\geq 5$ and
\begin{equation}
	\label{Totaro-1}
	\MSU_*[1/2]/I\simeq\mathbb{Z}[1/2][x_2,x_3,x_4].
\end{equation}

In \cite{H} by using of characteristic classes, Hoehn constructed a base sequence $W_1, W_2, W_3, \cdots $ of the rational cobordism ring $\MU^* \otimes \mathbb{Q}$ on which $\phi_{KH}$ has the values $A, B, C, D$, and $0$ for $W_i$ with $i>4$.

As another generalization of Oshanin's elliptic genus Schreieder in \cite{SC} studied a genus $\psi$ with logarithmic series  
$$\log_{\psi}(x)=\int_0^x\frac{dt}{R(t)},\,\,\,R(t)=\sqrt{1+q_1t+q_2t^2+q_3t^3+q_4t^4}.$$

The genus $\psi$ is easily calculable on  cobordism classes of complex projective spaces $\CP_i$, the generators of the domain $$\MU^*\otimes \mathbb{Q}=\mathbb{Q}[\CP_1,\CP_2, \cdots ].$$
This is because of the equation $(\log_{\psi}(x)')^2=1/R^2(x)=\sum_{i\geq 1}\psi(\CP_i)x^i$ we need only the Taylor expansion of $(1 + y)^{-1/2}$.

It is natural to ask whether one can calculate $\phi_{KH}$ in an elementary manner,
different from that relying on the formulas in \cite{H} and \cite{BU-PA}.

Viewed as a classifying map $\psi$ is strongly isomorphic to genus $\phi$ by the series $\mu(x)=\sum_{i\geq 0} \CP_ix^{i+1}$ \cite{B2}. 
This gives a method for explicit calculation  of $\phi_{KH}$.

In \cite{B, B1} we introduced the formal power series  
\begin{equation}
\label{A_{ij}}
A(x,y)=\sum A_{ij}x^iy^j=F(x,y)(x\omega(y)-y\omega(x)) \in \MU^*[[x,y]],
\end{equation} 
where 
$$F=F(x,y)=\sum\alpha_{ij}x^iy^j$$
is the universal formal group law over complex cobordism ring $\MU^*$ and 
$$
\omega(x)=\frac{\partial F(x,y)}{\partial y}(x,0)=1+\sum_{i\geq 1}w_ix^i
$$
is the invariant differential of $F$.
\medskip

The series $A(x,y)$ has proven to be interesting for the following reasons.

\begin{proposition}
	\label{prop.1}
	\cite{B1}. i) The obvious quotient map 
	$$f_B:\MU^*\to \MU^*/(A_{ij},\,i,j\geq 3)$$
	classifies a formal group law which is identical to the universal
	Buchstaber formal group law $F_B$, the universal formal group law of the form
	$$
	\frac{x^2A(y)-y^2A(x)}{xB(y)-yB(x)},
	$$ 
	where $A(0)=B(0)=1$.
	
	(ii) If $A'(0)=B'(0)$ then $B(x)$ is identical to the image of $\omega(x)$ under the  classifying map $f_B$.
	
\end{proposition}

\begin{proposition} \cite{B2}. 
	\label{prop.2}
	After it is tensored with rationals the classifying map $f_B$ of the Buchstaber formal group law is identical to the Krichever-Hoehn complex elliptic genus
	$$
	\phi: \MU^*\otimes \mathbb{Q}\to \Lambda_B \otimes \mathbb{Q}=\mathbb{Q}[\CP_1,\CP_2,\CP_3,\CP_4],
	$$
	where $\CP_i$ are complex projective spaces.	   
\end{proposition}

For explicit calculation of the Krichever-Hoehn genus the following observation is helpful.

\begin{proposition} \cite{B2}. 
	\label{prop.3}
	Over the ring $\MU^*\otimes \mathbb{Q}$ the series $\frac{x}{\omega(x)}$ is the strong isomorphism from the formal group law with logarithm series $$\int_0^x\frac{dt}{\sqrt{1+p_1t+p_2t^2+p_3t^3+p_4}}$$
	in \cite{SC} to the formal group law classified by $f$.
\end{proposition}

\bigskip

\bigskip

\medskip

\section{Some auxiliary combinatorial definitions}
\label{3}

By Euclid's algorithm for the natural numbers $m_1,m_2, \cdots, m_k$  one can find integers $\lambda_1, \lambda_2, \cdots \lambda_k $ such that 
\begin{equation}
	\label{lambda}
	\lambda_1m_1+\lambda_2m_2+\cdots  +\lambda_km_k=gcd(m_1,m_2,\cdots ,m_k).
\end{equation}

Let 
\begin{equation}
	\label{d}
	d(m)=gcd \bigg \{ \binom{m+1}{1},\binom{m+1}{2},\cdots ,\binom{m+1}{m-1}   \big|\,\, m\geq 1 \bigg \}.
\end{equation}

By \cite{KU} one has 
\begin{equation*}
	d(m)=
	\begin{cases} p, &\mbox{if $m+1=p^s$ for some prime  $p$,}\\
		1, &\mbox{otherwise}.
	\end{cases}
\end{equation*}

For the coefficients of universal formal group law $F(x,y)=\sum \alpha_{i\,j}x^iy^j$ the elements

\begin{equation}
\label{e_n}
e_{m}=\lambda_1\alpha_{1\,m}+\lambda_2\alpha_{2\,m-1}+\cdots +\lambda_m\alpha_{m\,1}	
\end{equation}  
are multiplicative generators in $\MU_*$.

By \cite{BU-U}, Theorem 9.9, or \cite{PANOV} 

\begin{align}
	\label{BU-U}
	&\frac{D(m)}{d(m)}=
	\begin{cases} d(m-1) &\mbox{if } m \neq 2^k-2,\\
		2&  \mbox{if } m=2^k-2.  
	\end{cases}	
\end{align}
where
\begin{equation}
	\label{D(m)}
	D(m)=
	gcd \bigg \{ \binom{m+1}{i}-\binom{m+1}{i-1}  
	\big|\,\, 2<i\leq m-1 , m\geq 5 \bigg \}.
\end{equation}

Let $m\geq 4$ and let $\lambda_2,\cdots ,\lambda_{m-2}$ are such integers that

\begin{equation}
	\label{d_2}
	d_2(m):=\sum_{i=2}^{m-2}\lambda_i\binom{m+1}{i}=gcd \bigg \{ \binom{m+1}{2},\cdots ,\binom{m+1}{m-2}  \bigg \}.
\end{equation}

Then by \cite{BU-U} Lemma 9.7 one has for $m\geq 3$

\begin{align} 
	\label{d_2d_1} 
	&d_2(m)=d(m)d(m-1).	
\end{align}

Note $d(n)$ are the Chern numbers of the generators in complex cobordism of dimension $2n$.

For the generators $$A_{ij}, \,\,i,j\geq 3,\,\, i+j-2=m$$ of the quotient ideal corresponding to $\Lambda_B$, the scalar ring of the universal Buchstaber formal group law in Proposition \ref{prop.1} and  the integers $\lambda_3,\cdots, \lambda_{m-1}$ corresponding to \eqref{D(m)} consider the linear combinations 

\begin{equation}
\label{T_n}	
T_{m}=\lambda_3 A_{3\,m-1}+\lambda_4 A_{4\,m-2}+\cdots \lambda_{m-1}A_{m-1\,3}.	
\end{equation}	

The elements $T_{p^s}$, where $p$ is a prime number,  and $e_i$ in \eqref{e_n} for $i\neq p^s$ will play a major role in Section \ref{realization}.

\medskip

\section{Realization of the universal Buchstaber formal group law localized away from 2.}
\label{realization}

Let $e_i$ and $T_i$ be as in \eqref{e_n} and \eqref{T_n} respectively.

Let $J_B$ be the ideal of $\MU_*=\mathbb{Z}[e_1,e_2, \cdots]$
$$J_B=\{A_{ij}, i,j\geq 3\},$$ 
the quotient ideal of the universal Buchstaber formal group law $F_B$ classified by	
$$f_B:\MU_*\to\MU_*/J_B=\Lambda_B.$$
The ideal $J_B$ is not prime as the quotient ring $\Lambda_B$
is not an integral domain: it has $2$-torsion element of  degree $12$ \cite{B1}. Here we use the results of \cite{BU-U}, that 
$\Lambda_{B}$ is generated by $f_{B}(e_j)$, $j=1,2,3,4$ and $j=p^r$, $r\geq 1$, $p$ is prime, and $j=2^k-2$, $k\geq 3$.
Then the ideal $Tor(\Lambda_B)$ is generated by the elements of order 2, namely $f_B(e_{j}),\,\,j=2^k-2$, $k \geq 3$.  

The ideal $Tor(\Lambda_B)$  is prime as $\Lambda_{\mathcal{B}}$ in
\begin{equation}
\label{integral domain}
f_{\mathcal{B}}:\MU_*\to \Lambda_B\to\Lambda_B/Tor(\Lambda_B)=\Lambda_{\mathcal{B}}.
\end{equation}
is
an integral domain and so is $J_{\mathcal{B}}=f_B^{-1}(Tor(\Lambda_B))$, the preimage ideal in $\MU_*$. Then the ideal
\begin{equation}
	\label{J_B}	
	J_{\mathcal{B}}=J_B+(e_{2^k-2}), k\geq 3
\end{equation}
is the kernel of the composition \eqref{integral domain}.
Denote by $F_{\mathcal{B}}$ the formal group law classified by $f_{\mathcal{B}}$. 

\bigskip

Let
\begin{align}
	\label{J}	
	J=\big(\mathcal{T}_l, l\geq 5 \big),  \text{ where } \mathcal{T}_l=
	\begin{cases}
		T_l, &\mbox{$n=p^s$, $p$ is a prime}  \\
		e_l, &\mbox{otherwise}.
	\end{cases}
\end{align}

Let $J_{\mathcal{B}}(l)\subset J_{\mathcal{B}}$, generated by those  elements whose degree is greater or equal $-2l$, and let $J(l)\subset  J $ be generated by $\mathcal{T}_5,\cdots, \mathcal{T}_l$.

\begin{remark}
	\label{remark-l}
	We note that those polynomials in $\mathbb{Z}[e_1,e_2,\cdots ,e_l]$ that are in the kernel of $f_{\mathcal{B}}$ can be viewed as the elements of $J_{\mathcal{B}}(l)$.
\end{remark}

\begin{proposition}
	\label{proposition-J_l}	 
	$J(l)=J_{\mathcal{B}}(l)$ for any natural $l\geq 5$. 	
\end{proposition}  

Proof. It is clear that $J(l)\subset J_{\mathcal{B}}(l)$.

Let us prove $ J_{\mathcal{B}}(l)\subset J(l)$ by induction on $l$. It is obvious for $l= 5$ as $T_5=A_{34}$. 

To prove $J_{\mathcal{B}}(l)=J_{\mathcal{B}}(l-1)+T_l$
note   
$$s_{i+j-2}(A_{ij})=\binom{i+j-1}{j-1}-\binom{i+j-1}{j}.$$ 
Indeed, modulo decomposable elements  
$$A_{ij}=\alpha_{i-1j}-\alpha_{ij-1}$$ 
and $s_{i+j-1}(\alpha_{ij})=-\binom{i+j}{i}$. Now apply Euclid's algorithm for  $m_i=s_l(A_{i,l+2-i})$, fix the integers $\lambda_i$ and consider  the elements $T_l$ in \eqref{T_n}.

The combinatorial identities in Section 3 implies that 
\begin{equation}
\label{s_l}
s_l(T_l)=D(l)
\end{equation}
is
the greatest common divisor of the
integers $s_l( A_{ij})$ for $A_{ij}\in J_B$. $i+j-2=l$. 

It follows that 
$$A_{ij}=\frac{s_n(A_{ij})}{D(l)}T_l + P(e_1,e_2,\cdots, e_{l-1}),$$
for some polynomial $P$, i.e., 
$$A_{ij}=P(e_1,e_2,\cdots ,e_{l-1}) \text{ modulo } T_l $$ 
Therefore  $P(e_1,e_2,\cdots ,e_{l-1})$ is in the kernel of $f_{\mathcal{B}}$, i.e., is in $J_{\mathcal{B}}(l-1)=J(l-1)$ by above Remark \ref{remark-l}.
\qed

Let
$A_l=\mathbb{Z}[e_1,e_2,\cdots , e_l]$ and $A^{l+1}=\mathbb{Z}[e_{l+1}, e_{l+2}, \cdots]$ i.e., $\MU_*=A_l\otimes A^{l+1}$.   
Let $J(l)$ as above be generated by $\mathcal{T}_1,\cdots, \mathcal{T}_l$. The preimage of $J(l)$ by obvious inclusion defines the ideal of $A_l$ denoted by same symbol so that 
$$\MU_*/J(l)=A_l/J(l)\otimes A^{l+1}.$$

\medskip

\begin{proposition}
\label{regular}
	i) The ideal $J$ in \eqref{J} is regular;
	
	ii) $A_l/J(l)$ and $\MU_*/J(l)$ are integral domains, or equivalently $J(l)$ is prime.  
\end{proposition}

It is clear that ii) implies i): If $\MU_*/J(l)$ is integral domain for any $l\geq 5$, i.e., it has no zero divisors,  then multiplication by $\mathcal{T}_{l+1}$ is monorphism. Therefore the sequence $\mathcal{T}_5,\mathcal{T}_6, \cdots $ of generators of $J$ is regular.

We will see that ii) follows from the proof of Proposition 6.5 in \cite{BU-U} and the following

\begin{lemma} For  $p^r\leq l<p^{r+1}$ the ring $A_l/J(l)\otimes \mathbb{F}_p$ is additively generated by the following monomials 
	\label{monomials}
	
	For $p=2$,
	$$
	\alpha_1^{m_1}\beta_2^{m_2}\alpha_3^{m_3}\beta_2^{k_1}\beta_{2^2}^{k_2}\cdots \beta_{2^{r-1}}^{k_{r-1}}\beta_{2^r}^m, \,\,\,\,\,\,\,
	k_1,k_2,k_{r-1}=0,1;
	$$
	
	\medskip
	
	For $p=3$,
	$$
	\alpha_1^{m_1}\beta_2^{m_2}\beta_4^{m_4}\beta_3^{k_1}\beta_{3^2}^{k_2}\cdots \beta_{3^{r-1}}^{k_{r-1}}\beta_{3^r}^m, \,\,\,\,\,\,\,
	k_1,k_2,k_{r-1}=0,1,2;
	$$
	
	\medskip
	
	For prime $p>3$,
	$$
	\alpha_1^{m_1}\beta_2^{m_2}\alpha_3^{m_3}\beta_4^{m_4}\beta_p^{k_1}\beta_{p^2}^{k_2}\cdots \beta_{p^{r-1}}^{k_{r-1}}\beta_{p^r}^m, \,\,\,\,\,\,\,
	0\leq k_1,k_2,k_{r-1}\leq p-1,
	$$
	not divisible by $\alpha_1^{j_1}\beta_2^{j_2}\alpha_3^{j_3}\beta_4^{j_4}$, where $(j_1j_2j_3j_4)$ corresponds to leading lexicographical monomial ordering for which 
	$\lambda_{j_1j_2j_3j_4}\not\equiv 0 \mod p$ in 
	$$
	p\beta_p=\sum \lambda_{j_1j_2j_3j_4}\alpha_1^{j_1}\beta_2^{j_2}\alpha_3^{j_3}\beta_4^{j_4}.
	$$
\end{lemma}  

Proof. We follow the proof of Proposition 6.5 in \cite{BU-U}.
To get the generating monomials in Lemma \ref{monomials} we need only to modify the generating monomials of $\Lambda_{\mathcal{B}}\otimes \mathbb{F}_p$. In particular,  in (6.17), (6.18) and (6.20) there are extra factors for $A_l/J(l)$, $p^r\leq l<p^{r+1}$, namely
$$\beta_{p^{r+1}}^{k_1}\beta_{p^{r+2}}^{k_2}\cdots \beta_{p^{r+s}}^{k_s}\cdots, \,\,\,\,  k_1,k_2,\cdots \leq p-1.$$
We have to replace these factors by 
$$\beta_{p^r}^{pk_1+p^2k_2+\cdots +p^sk_s+\cdots }$$
This is because of $\beta_{p^{r+1}}\notin A_{l}$. 
In this way for each $m$ we keep the total number of generating monomials of $\Lambda^{-2m}\otimes \mathbb{F}_p$ since there is no relation $(6.2)$ of \cite{BU-U} in our ring $A_l/J(l)\otimes \mathbb{F}_p$. 
\qed

\bigskip

Denote by $[\mathcal{T}_i]$ the cobordism class representing the generator $\mathcal{T}_i$ of the ideal $J$. Consider the sequence 
$$\Sigma_T=([\mathcal{T}_5], [\mathcal{T}_6], \cdots).$$

The Sullivan-Baas construction \cite{BA}  of cobordism with singularities 
$\Sigma_{T}$ gives a cohomology theory $\MU^*_{\Sigma_T}(-)$ which by regularity 
of the ideal $J$ has a scalar ring 
$$\MU^*_{\Sigma_T}(pt)=\MU_*/J=\Lambda_{\mathcal{B}}.$$

By Mironov \cite{MI} (Theorem 4.3 and Theorem 4.5) $\MU^*_{\Sigma_T}(-)$ admits an associate multiplication and all obstructions to commutativity are in $\Lambda_\mathcal{B}\otimes \mathbb{F}_2$. Therefore after localization away from 2 
all obstructions vanish and we get a commutative cohomology $$\emph{h}_{\mathcal{B}}^*(-):=\MU^*_{\Sigma_T}[1/2](-).$$ 
Here we recall that $\Lambda_{B}[1/2]=\Lambda_{\mathcal{B}}[1/2]$ by definition of $\Lambda_{\mathcal{B}}$.

It is clear that $\emph{h}_{\mathcal{B}}^*(-)$ is complex oriented as the  Atiyah-Hirzebruch
spectral sequence 
$H^*(-,\emph{h}_{\mathcal{B}}^*(pt))\Rightarrow \emph{h}_{\mathcal{B}}^*(-)$
collapses for $\BU(1)\times \BU(1)$. 

\bigskip

Thus we can state

\begin{theorem}
\label{main-lemma2}	
There exist a commutative complex oriented cohomology 
$h_{\mathcal{B}}^*(-)$ with scalar ring isomorphic to  $\Lambda_{\mathcal{B}}[1/2]$, the ring of coefficients of the universal Buchstaber formal group law localized away from 2.  
\end{theorem}

This result without 
a complete proof, is announced in short communications of MMS \cite{B-4}.

\bigskip

 \section{The restriction of the Buchstaber genus on $\MSU_*[1/2]$}
 \label{5}

Taking into account \cite{NOV} that after localized away from 2, the forgetful map from the special unitary cobordism $\MSU_*[1/2]=\mathbb{Z}[1/2][x_2,x_3,\cdots]$ to complex cobordism $\MU_*[1/2]$ is an injection, define the following  ideal extensions in $\MU_*][1/2]$:

$J_{SU}^e$,  generated by any polynomial generators $x_n$ of $\MSU_*[1/2]$,   $n\geq 5$ viewed as elements in $\MU_*[1/2]$ by forgetful injection map;

$J_{T}^e$,  generated by $SU$-flops \cite{TO} of dimension $\geq 10$ again viewed as elements in $\MU_*[1/2]$;

$J_B^e$,  the contraction ideal by the obvious inclusion $\MU_*\to \MU_*[1/2]$ of the ideal $J_B$ of $\MU_*$ generated by the elements $\{A_{ij}, i,j\geq 3\}$, defined in Section 3.

\begin{proposition}
	\label{main-lemma}
	i) $J^e_B=J_{T}^e=J_{SU}^e$; ii)
	When restricted on $\MSU_*[1/2]$ the classifying map of the Buchstaber formal group law localized away from 2, gives a genus with the scalar ring $\mathbb{Z}[1/2][x_2,x_3,x_4]$, $|x_i|=2i$.
\end{proposition}
One motivation is the restricted Krichever-Hoehn complex elliptic genus below, studied in \cite{H} and \cite{TO}. Another construction with the scalar ring $\mathbb{Z}_{(2)}[a,b]$, $|a|=2$, $|b|=6$ see in \cite{BU-BU}. 

Proof.
$J_{T}^e=J_{SU}^c\subseteq J^e_B$:  The homomorphism 
$$\MU_* \xrightarrow{f_B} \Lambda_B \xrightarrow{\subset} \Lambda_B\otimes \mathbb{Q}= \mathbb{Q}[x_1,x_2,x_3,x_4]
$$
is a  specialization of the complex elliptic genus
$$
\MU_* \xrightarrow{\subset} \MU_*\otimes \mathbb{Q}\to \mathbb{Q}[x_1,x_2,x_3,x_4]
$$
by \cite{B}, therefore vanishes on the kernel of the complex elliptic genus which is the ideal $I$ of complex flops by \cite{TO} (Theorem 4.1). 
On the other hand  the ring $\Lambda_B[1/2]$  is torsion free and injected in $\mathbb{Q}[x_1,x_2,x_3,x_4]$ by \cite{BU-U}. So the ring homomorphism 
$\MU_*[1/2] \xrightarrow {f_{B}} \Lambda_{B}[1/2]$ 
vanishes on $I$, therefore it vanishes on $SU$-flops. Moreover by \cite{TO} the ideal of $SU$-flops $J_{T}$ in $\MSU_*[1/2]$ contains the polynomial generators $x_n,\,\,\,n\geq 5$ constructed by using Euclid's algorithm and $SU$-flops. Therefore $J_{T}=(x_5,\cdots )$ and $J_{T}^e=J_{SU}^e\subseteq J_B$. 

To prove $J^e_B\subseteq J_{SU}^e$ note by \eqref{s_l} 
 $T_n$ satisfies the criteria for the membership of the set of polynomial generators in  
$$\MSU^*[1/2]=\mathbb{Z}[1/2][x_2,x_3,\cdots ],\,\,|x_i|=2i,$$ 
described by Novikov in \cite{NOV}.
In particular, 
an $SU$-manifold $M$ of real dimension $2n$, $n\geq 2$ is a polynomial
generator if and only if $s_n(M)$, the main Chern characteristic number is as follows
\begin{align}
	\label{Novikov}
	s_n(M)=
	\begin{cases} \pm 2^kp &\mbox{if } n=p^l,\,\,\,p\,\,\, \mbox{is odd prime, }\\
		\pm 2^kp &\mbox{if } n+1=p^l,\,\,\, p\,\,\,\mbox{is odd prime, } \\
		\pm 2^k& \mbox{otherwise}.
	\end{cases}
\end{align}

 It follows that the generators $A_{ij}$ of  $J_B$, with $i+j=n+2$ and the generators $x_n$ of $J_{SU}^e\subset J^e_B$ are related as follows
$A_{ij}\equiv \frac{s_n(A_{ij})}{D(n)} T_n$, $T_n\equiv \pm 2^{k(n)}x_n$,  $\mod$ decomposables. 
Therefore we can proceed as in the proof of Proposition \ref{proposition-J_l}.  
\qed

\bigskip

 Recall also that $\Lambda_B[1/2]=\Lambda_{\mathcal{B}}[1/2]$ is an integral domain. Note that Proposition \ref{regular} implies 

\begin{proposition}
	\label{regular-B}
	The sequence $\{x_n\}, n\geq 5$ of any polynomial generators in $\MSU_*[1/2]$ viewed as elements in $\MU_*[1/2]$ by forgetful map is regular. 
\end{proposition}

\begin{corollary}
	\label{cor-1} After restriction on $\MSU_*[1/2]$ the Buchsteber genus $f_{\mathcal{B}}$ gives a cohomology theory $\MSU^*_\Sigma [1/2]$  with singularities $\Sigma=(x_5, \cdots )$, with the scalar ring 
	$$
	\mathbb{Z}[1/2][x_2,x_3,x_4], \,\,\,|x_i|=2i.
	$$
\end{corollary}

\bigskip

\section{$\MSU_*[1/2]$}
\label{6}

Let  $F_U(x,y)=\sum\alpha_{ij}u^iv^j$ be the universal formal group law. Recall the idempotent in \cite{BU-1}, \cite{C-P} 

$$\pi_0:\MU_*\to \MU_* : \,\,\,  \pi_0=1 + \sum_{k \ge 2} \alpha_{1k} \partial_k$$
and the projection $\pi_0:\MU_*\to W_*=Im \pi_0$.  

Then $W_*$ is a ring with multiplication $*$ and  $\pi_0(a\cdot b)=a*b$. By \cite{C-P} Proposition 2.15  the multiplication $*$ is given by
$$
a*b = ab + 2[V]\partial a \partial b, 
$$
where $[V]=\alpha_{12}\in\MU_4$ is the cobordism class $\CP_1^2 - \CP_2$.

$W_*$ is complex oriented and by \cite{C-P} Proposition 3.12 the ring $W_*[1/2]$ is generated by the coefficients of the formal group law 
$$
F_W=F_W(x,y)=\sum w_{ij}x^iy^j.$$  
 
Following \cite{BU-1}, \cite{C-P} one can calculate $w_{ij}$ in terms of $\alpha_{kl}$ as follows. Consider the multiplicative cohomology theory $\Gamma $ with
$$\pi_*(\Gamma)=\MU_*[t]/(t^2-\alpha_{11}t-2\alpha_{21}),$$
the free $\MU_*$ module generated by $1$ and $t$.

There is a natural multiplicative transformation $\phi: W \to \Gamma$ given by $$\phi_*(x) = x + t\partial x$$
for $x\in W_*$. The restriction of  $\phi_*[1/2]$ on $\MSU_*[1/2]$, the subring in $W_*[1/2]$ of cycles of $\partial $,  is the natural inclusion in $\MU_*[1/2]$.  

Then  
\begin{equation*}
\phi_*F_W=u+v+\sum_{i,j\geq 1}(w_{ij}+t\partial w_{ij})u^iv^j
\end{equation*}
is strongly isomorphic to $F_{U}$ (considered as a formal group low over $\pi_*(\Gamma)$ via the natural inclusion) by a series $\gamma^{-1}$, i.e.,
\begin{equation}
\label{calc-2}
u+v+\sum_{i,j\geq 1}(w_{ij}+t\partial w_{ij})u^iv^j=\gamma F_U(\gamma^{-1}(u),\gamma^{-1}(v)).
\end{equation}

Finally, we need to apply for $\gamma$ in \eqref{calc-2} Lemma 3 in \cite{BU-1}, which says that
any orientation $w\in W_2$ gives the following identity in $\pi_*(\Gamma)$ 
\begin{equation}
\label{gamma}
\gamma(u)=\phi_*(w)=\pi_0(u)+t\partial u=u+tu\bar{u}+
\sum_{i\geq 2} \alpha_{i1}u\bar{u}^i.
\end{equation}

\medskip

By \cite{C-P} one can specify an orientation of $W$ such that 
$$gcd(w_{ij}, \,i+j-1=k) =d(k)d(k-1)$$ modulo a power of 2. This allows to construct the generators of $W_*[1/2]$.

\medskip

In particular, 
one can calculate main Chern numbers $s_k(w_{ij})$ in terms of main Chern numbers of $\alpha_{ij}$ as follows. By \cite{C-P} Lemma 3.5
any orientation $w\in W_2$ gives $w_i\in W_{2i}$ and $\lambda\in \MU_2=W_2$ such that one has in $\pi_*(\Gamma)$ 
\begin{equation*}
	\gamma(u)=u-(\lambda+(2l+1)t)u^2+\sum_{i\geq 2}\gamma_{i+1}u^{i+1} mod J^2+tJ,
\end{equation*}
where $2l=\partial\lambda$, $l\in \mathbb{Z}$, $\gamma_{i+1}=(-1)^i\alpha_{1i}+w_i$,
$J$ is the ideal in $\MU_*$ of elements of positive degree.

\medskip

\begin{lemma} 
	\label{lem}
Let $k=i+j-1\geq 3$  is not of the form $k=2^l=p^s-1$ for some odd prime $p$.  There is a choice of complex orientation $w$ for the theory $W$ such that  
	\begin{align*}
		s_k(w_{ij})=
		\begin{cases} p s_k(\alpha_{ij}) &\mbox{if } k=p^s,\,\,\,p\,\,\, \mbox{is any prime, $s>0$},\\
			s_k(\alpha_{ij}) &\mbox{otherwise}.\\
		\end{cases}	
	\end{align*}
\end{lemma}

Proof. By using \eqref{gamma} it is proved in \cite{C-P} (Lemma 3.9) that for $k\geq 3$ one has  modulo decomposable elements

\begin{align} 
	\label{f-1}
	&\phi_*(w_{1k})=\alpha_{1k}+(k+1)((-1)^k\alpha_{1k}+w_k);\\
	\label{f-2}
	&\phi_*(w_{ij})=\alpha_{ij}+(-1)^k\binom{k+1}{i}\alpha_{1k}+\binom{k+1}{i}w_k;
\end{align}

Choose a  complex orientation $w$ for the theory $W$ such that  elements $w_k$ satisfy the following conditions

\begin{equation}
\label{f-2}
1+(-1)^k(k+1)-s_k(w_k)=d(k-1)
\end{equation}

The it is easily checked that $w_k$ indeed belongs to $W_*$, that is $s_k(w_k)$ is divisible by $d(k-1)d(k)$.

Then by \eqref{f-1} and \eqref{f-2} we have for $k=p^s$

\begin{align*}
	w_{1k}=&-k\alpha_{1k}+(k+1)w_k=-k\alpha_{1k}+\alpha_{1k}(p^s+p)=p\alpha_{1k};\\
	w_{ij}= &\alpha_{ij}-\binom{k+1}{i}\alpha_{1k}+\binom{k+1}{i}w_k\\
	=&\alpha_{ij}-\alpha_{ij}(k+1)+\alpha_{ij}(k+p)=p \alpha_{ij}.
\end{align*}

\medskip

For $k=2^l$ one has

\begin{align*}
w_{1k}=&(2+k)\alpha_{1k}+(k+1)w_k=(2+k)\alpha_{1k}-\alpha_{1k}k=2\alpha_{1k};\\
w_{ij}=&\alpha_{ij}+\binom{k+1}{i}\alpha_{1k}-\binom{k+1}{i}w_k
=\alpha_{ij}+\alpha_{ij}(k+1)-k \alpha_{ij}=2\alpha_{ij}.
\end{align*}

\medskip

Similarly for other cases.

\qed

\medskip

If $k=2^l=p^s-1\geq 3$ for some odd prime $p$, then by \cite{C-P} Lemma 3.15 one has $k=8$ or $k=2^{2^n}$.  We have to replace $d(k-1)=2$ in \eqref{f-2} by $4$ , if $k=8$ to get $\phi_*(w_{i\,9-i})=4\alpha_{i\,9-i}$ and by $-2^{2^n}$, if $k=2^{2^n}$ to get $\phi_*(w_{i\,2^{2^n}+1-i})=-2^{2^n}\alpha_{i\,2^{2^n}+1-i}$.

\medskip

Together with Lemma \ref{lem}, this implies

\begin{corollary}
\label{b_k}
 Let $k\geq 3$. By \eqref{lambda} and \eqref{d} let $\lambda_i$ be such that the linear combination ${\bf a}_k=\sum_{i=1}^{k}\lambda_i\alpha_{ik+1-i}$ is a polynomial generator in $\MU_{2k}$. Then
$$
{\bf b}_k=\sum_{i=1}^{k}\lambda_iw_{ik+1-i}.
$$
is a polynomial generator in $W_{2k}[1/2]$.
\end{corollary}

\bigskip

Then  $\MSU_*[1/2]$ is a subring of cycles of the boundary operation $\partial $ in $W_*[1/2]$ with multiplication $*$
$$\partial (a*b)=a*\partial b+\partial a *b-\CP_1\partial a*\partial b.$$

One has $\partial \CP_1=2$ and $a*b=a\cdot b $ whenever $a\in Im \partial $ or  $b\in Im \partial $. Therefore 
$$
\partial(\CP_1 * \alpha)=2\alpha-\CP_1\cdot \partial \alpha, \forall \alpha \in W. 
$$
As mentioned in \cite{L-C-P} this implies that 
\begin{equation}
\label{su-component}
\alpha=1/2\partial(\CP_1 * \alpha)+1/2\CP_1\cdot \partial \alpha,
\end{equation}
and therefore  $W[1/2]$ is generated by $1$ and $\CP_1$  as a $\MSU_*[1/2]$ module . It is easily checked that this module is free.

\begin{proposition} 
\label{msu-b_k}
Let ${\bf b}_k$ be as in Corollary \ref{b_k}. Then
$$\MSU_*[1/2]=\mathbb{Z}[1/2][x_2, x_k:k\geq 3],$$
where $x_2=\CP_2-9/8\CP_1^2$ and  $x_k=\partial(\CP_1 *{\bf b}_k)$.
\end{proposition}

Proof. One has for the values of the Chern numbers
	\begin{align*}
		&c_1c_1[\CP_1^2]=8, &&c_1c_1[\CP_2]=9, &&&c_2[\CP_1^2]=4, &&&&& c_2[\CP_2]=3.\\ 
	\end{align*}
	This imply that that $c_1c_1[{\bf b}_2]=0$. There are no more Chern numbers having $c_1$ as a factor and  $s_2[{\bf b}_2]=s_2[\CP^2]=3$.  Therefore ${\bf b}_2$ forms a generator of $\MSU_4[1/2]$.
	
	Apply \eqref{su-component}. The main Chern number vanishes on the  second (decomposable) component of  
	$$
	{\bf b}_k=1/2\partial(\CP_1 * {\bf b}_k)+1/2\CP_1\cdot \partial {\bf b}_k, 
	$$
i.e., the first component  $x_k$ has the  main Chern number $2s_k({\bf b}_k)$. 	
\qed

\bigskip

 \bigskip

 \section{The restriction of the classifying map of $F_{Ab}$ on $\MSU^*[1/2]$.}
 \label{7}
 
 As above let $F_U=\sum \alpha_{ij}x^iy^j$ be the universal formal group law. By definition the coefficient ring of the universal abelian formal group law $F_{Ab}$ is the quotient ring
 \begin{equation}
 	\label{Ab}
 	\Lambda_{Ab}=\MU_*/I_{Ab}, \text { where } I_{Ab}=(\alpha_{ij}, i,j >1).
 \end{equation}
 
 Let us apply Euclid's algorithm for the Chern numbers $s_{m-1}(\alpha_{i,m-i})$ in \eqref{d_2}
 
 Let 
 $$z_{k}=\sum_{i=2}^{k-1}\lambda_i\alpha_{i\, k+1-i},\,\,\,k\geq 3. $$   
 
 By \cite{B-KH}, \cite{BUSATO} one has $I_{Ab}=I_{AB}=(z_k, \,\,k\geq 3)$.
 

 Consider the composition
 \begin{equation}
 	r_{Ab}:\MSU_*[1/2] \xrightarrow{\subset} \MU_*[1/2] \xrightarrow{f_{Ab}} \Lambda_{Ab}[1/2],
 \end{equation}
 where $\subset$ is forgetful map.
 
 \begin{proposition}
 One has the following polynomial generators in $\MSU_*[1/2]$ viewed as the elements in $\MU_*[1/2]$ 
 \begin{align*}
 		&x_2=\CP_2-\frac{9}{8}\CP_1^2,	
 		&&x_3=-\alpha_{22},
 		&&x_4=-\alpha_{23}-\frac{3}{2}x_3\CP_1.\\
 \end{align*}

\end{proposition}
To prove this we have to check that all Chern numbers of $x_i$ having factor $c_1$ are zero. Then we have to check the main Chern number $s_i(x_i)$ for Novikov's criteria. 
 
We already did this for $x_2$ in the proof of Proposition \ref{msu-b_k}. Then by definition $x_3$ is the coefficient $-\alpha_{22}\in I_{AB}$ of the universal formal group law. In $\MU_*$ one has
 
 \begin{align*}
 	&2\alpha_{22}=-3\CP_3+8\CP_1\CP_2-5\CP_1^3,\\
 	&\alpha_{23}=2\CP_1^4-7\CP_1^2*\CP_2+3\CP_2^2+4\CP_1\CP_3-2\CP_4.
\end{align*}	

Let us compute the Chern numbers of $x_3=-\alpha_{22}$.  One has

 \begin{align*}
 	& X             &c_3(X)   && c_1c_2(X)     &&& c_1c_1c_1(X)\\
 	&\CP_3           &4        && 24            &&& 64  \\
 	&\CP_1\CP_2     &6        && 24            &&& 54    \\
 	&\CP_1^3        &8        && 24            &&& 48.\\
 \end{align*}

 	It follows all Chern numbers of $\alpha_{22}$ having factor $c_1$ are zero. 
 	Then $\alpha_{22}$ forms a generator in $\MSU_6[1/2]$ as $s_3[-2\alpha_{22}]=4\cdot 3$.
  
 Similarly for $x_4$:  the main Chern number $s_4(x_4)=2\cdot 5$ fits for Novikov's criteria and one has

 \begin{align*}
 	&X          &c_1c_1c_1c_1(X)  &&c_1c_1c_2(X) &&&c_1c_3(X)    \\
 	&\CP_1^4     &384             &&192          &&&64              \\
 	&\CP_1^2\CP_2 &432            && 204         &&&60           \\  
    & \CP_2^2     &486            && 216         &&&54           \\
    &\CP_1\CP_3   &512            && 224         &&&56            \\
    &\CP_4       &625             && 250         &&&50.              \\
\end{align*}

 Note $F_{Ab}$ is a specialization of the Buchstaber formal group law $F_B$.  In particular one can put $A(x)=B(x)^2$ in Proposition \eqref{prop.1} to specify $F_B$ to $F_{Ab}$
 over  torsion free ring $\Lambda_{Ab}[1/2]$. Then Proposition \ref{regular-B} implies
 
 \begin{proposition} 
 	After restriction on $\MSU_*[1/2]$ the classifying map of the universal abelian formal group law becomes the one-parameter genus
 	\label{r}
 	$$r_{Ab}: \MSU_*[1/2] \to \mathbb{Z}[1/2][x_2].$$	 
 \end{proposition}


\end{document}